\documentclass{amsart}

\usepackage{fancyhdr}
\usepackage{lastpage}
\usepackage{stmaryrd,yhmath}

\newcommand{\bdis}{\begin{displaymath}}
\newcommand{\edis}{\end{displaymath}}
\newcommand{\be}{\begin{equation}}
\newcommand{\ee}{\end{equation}}
\newcommand{\mbb}{\mathbb}
\newcommand{\mcal}{\mathcal}

\newcommand{\vp}{\varphi}

\newcommand{\vth}{\vartheta}

\newcommand{\mT}{\mathring{T}}

\newcommand{\zf}{\zeta\left(\frac{1}{2}+it\right)}

\theoremstyle{definition}

\theoremstyle{remark}
\newtheorem{remark}[]{Remark}

\newtheorem*{mydef1}{{\bf Theorem}}

\numberwithin{equation}{section}

\begin{document}

\title{Jacob's ladders and certain asymptotic multiplicative formula for the function $|\zf|^2$}

\author{Jan Moser}

\address{Department of Mathematical Analysis and Numerical Mathematics, Comenius University, Mlynska Dolina M105, 842 48 Bratislava, SLOVAKIA}

\email{jan.mozer@fmph.uniba.sk}


\begin{abstract}
In this paper it is proved that a mean-value of the product of some factors $|\zeta|^2$ is asymptotically equal to the product of the mean-values
of $\zeta|^2$, and this holds true for every fixed number of the factors.
\end{abstract}

\maketitle

\section{Introduction}

\subsection{}

Let

\be \label{1.1}
Z(t)=e^{i\vth(t)}\zf,\ \vth(t)=-\frac t2\ln\pi+\text{Im}\ln\Gamma\left(\frac 14+i\frac t2\right)
\ee
be the signal generated by the Riemann zeta-function. Hardy and Littlewood started to study the following integral in 1918

\be \label{1.2}
\int_0^T\left|\zf\right|^2{\rm d}t=\int_0^T Z^2(t){\rm d}t ,
\ee

and they have derived the following formula (see \cite{1}, pp. 122, 151-156)

\be \label{1.3}
\int_0^T Z^2(t){\rm d}t\sim T\ln T,\quad T\to\infty .
\ee

In this direction the Titchmarsh-Kober-Atkinson (TKA) formula

\be \label{1.4}
\int_0^\infty Z^2(t)e^{-2\delta t}{\rm d}t=\frac{c-\ln(4\pi\delta)}{2\sin\delta}+\sum_{n=0}^N c_n\delta^n+\mcal{O}(\delta^{N+\epsilon})
\ee

(see \cite{5}, p. 141) where $c$ is the Euler's constant, remained as an isolated result for the period of 56 years. In our paper \cite{2}, we have
discovered the nonlinear integral equation

\be \label{1.5}
\int_0^{\mu[x(T)]} Z^2(t)e^{-\frac{2}{x(T)}}{\rm d}t=\int_0^T Z^2(t){\rm d}t ,
\ee
where
\bdis
\mu(y)\geq 7y\ln y,\ \mu(y)\to y=\vp_\mu(T)=\vp(T) ,
\edis
in which the essence of the TKA formula (\ref{1.4}) is inclosed. Namely, we have shown in \cite{2} that the following infinite set of the almost-exact
expressions of the Hardy-Littlewood integral (\ref{1.1}) takes place

\be \label{1.6}
\begin{split}
& \int_0^T Z^2(t){\rm d}t=\vp_1(T)\ln\vp_1(T)+(c-\ln 2\pi)\vp_1(T)+c_0+\mcal{O}\left(\frac{\ln T}{T}\right), \\
& \vp_1(T)=\frac 12\vp(T)
\end{split}
\ee
where $\vp(T)$ (and, of course, also $\vp_1(T)$) is the Jacob's ladder, i.e. an arbitrary solution to the nonlinear integral equation (\ref{1.5}).

\begin{remark}
We have proved in the paper of reference \cite{2} that except the asymptotic formula (\ref{1.3}) possessing an unbounded error term (let us remark
in this direction that the formulae of Ingham, Titchmarsh and Balasubramanian also posses the unbounded errors, comp. \cite{2}, Remark 2) there is
an infinite family of almost exact representations (\ref{1.6}) of the Hardy-Littlewood integral (\ref{1.3}).
\end{remark}

\begin{remark}
We can formulate the result (\ref{1.6}) as follows: the Jacob's ladders $\vp_1(t)$ are the asymptotic solutions of the transcendental equation
\bdis
\int_0^T Z^2(t){\rm d}t=V(T)\ln V(T)+(c-\ln 2\pi)V(T)+c_0 .
\edis
\end{remark}

\subsection{}

Next, we have proved (see \cite{3}, (8.3)) the following sixth-order formula

\be \label{1.7}
\begin{split}
& \int_T^{T+U_1}\left|\zeta\left(\frac 12+i\vp_1(t)\right)\right|^4\left|\zf\right|^2{\rm d}t\sim \frac{1}{2\pi^2}U_1\ln^5T, \\
& U_1=T^{\frac 78+2\epsilon} , \quad T\to\infty ,
\end{split}
\ee
and (see \cite{4}, (2.5)) the following eight-order formula

\be \label{1.8}
\begin{split}
& \int_{\mT}^{\widering{T+U_2}}\left|\zeta\left(\frac 12+i\vp_2(t)\right)\right|^4\left|\zf\right|^4{\rm d}t\sim\frac{1}{4\pi^4}U_2\ln^8T, \ T\to\infty
\end{split}
\ee
where
\bdis
[T,T+U_2]=\vp_2\left\{\left[ \mT,\widering{T+U_2}\right]\right\},\quad U_2=T^{\frac{13}{14}+2\epsilon} ,
\edis
and $\vp_2$ is the Jacob's ladder of the second order.

\subsection{}

A motivation for the next step is in the well-known multiplicative formula

\be\label{1.9}
M\left(\prod_{k=1}^n X_k\right)=\prod_{k=1}^n M(X_k)
\ee
from the theory of probability, where $X_k$ are independent random variables and $M$ is the population mean. \\

We shall introduce, in connection with the formula (\ref{1.9}), some new classes of non-local formulae for the function $|\zf|^2$. Namely, we shall show
that there are functions
\bdis
w_k(t),\quad k=0,1,\dots ,n,\qquad t\in [T,T+U]
\edis
such that the following formula
\be \label{1.10}
\begin{split}
& \frac 1U\int_T^{T+U}\prod_{k=0}^n \left|\zeta\left(\frac 12+iw_k(t)\right)\right|^2{\rm d}t\sim \\
& \sim\prod_{k=0}^n\frac{1}{w_k(T+U)-w_k(T)}\int_{w_k(T)}^{w_k(T+U)}\left|\zf\right|^2{\rm d}t,\quad T\to\infty
\end{split}
\ee
holds true for every fixed $n\in \mbb{N}$.

\section{The result}

Let

\be \label{2.1}
\begin{split}
& y=\frac 12\vp(t)=\vp_1(t);\quad \vp_1^0(t)=t,\ \vp_1^1(t)=\vp_1(t), \\
& \vp_1^2(t)=\vp_1(\vp_1(t)), \dots , \vp_1^k(t)=\vp_1(\vp_1(\dots(\vp_1(t)))\dots),\ \dots , \\
& t\in [T,T+U] ,
\end{split}
\ee
where $\vp_1^k(t)$ denotes the $k$-th iteration of the Jacob's ladder
\bdis
y=\vp_1(t),\quad t\geq T_0[\vp_1] .
\edis
The following Theorem holds true.

\begin{mydef1}
Let
\be \label{2.2}
U\in\left(\left. 0,\frac{T}{\ln^2T}\right]\right. .
\ee
Then for every fixed $n\in\mbb{N}$ we have
\be \label{2.3}
\begin{split}
& \frac 1U\int_T^{T+U}\prod_{k=0}^n\left|\zeta\left(\frac 12 +i\vp_1^k(t)\right)\right|^2{\rm d}t\sim \\
& \sim\prod_{k=0}^n\frac{1}{\vp_1^k(T+U)-\vp_1^k(T)}\int_{\vp_1^k(T)}^{\vp_1^k(T+U)}\left|\zf\right|^2{\rm d}t,\quad T\to\infty
\end{split}
\ee
where
\be \label{2.4}
\vp_1^k(t)\geq (1-\epsilon)T,\ k=0,1,\dots ,n+1,\ t\in [T,T+U] .
\ee
\end{mydef1}

\begin{remark}
Some nonlocal and non-linear interaction of the signals
\bdis
\left|\zf\right|^2,\ \left|\zeta\left(\frac 12 +i\vp_1^1(t)\right)\right|^2,\ \dots \ , \left|\zeta\left(\frac 12 +i\vp_1^{n}(t)\right)\right|^2
\edis
is expressed by the formula (\ref{2.3}) and, simultaneously, a new art of the asymptotic independence of the partial functions
\bdis
\left|\zf\right|^2,\quad t\in [\vp_1^k(T),\vp_1^k(T+U)],\ k=0,1,\dots \ .
\edis
is expressed by this formula.
\end{remark}

\begin{remark}
We can put in the formula (\ref{2.3}) for example
\bdis
U=\frac{1}{T^2},\ U=\frac{2\pi}{\ln T}\sim t_{\nu+1}-t_\nu,\ t_\nu\in\left[ T,\frac{T}{\ln^2T}\right],\ \dots
\edis
where $\{ t_\nu\}$ is the Gram sequence, i.e. we can begin from the microscopic segments
\bdis
[T,T+U]=\left[ T,T+\frac{1}{T^2}\right],\ \left[ T,T+\frac{2\pi}{\ln T}\right],\dots \ .
\edis
\end{remark}

\begin{remark}
It is obvious that the formula (\ref{2.3}) cannot be obtained in the theories of Balasubramanian, Heath-Brown and Ivic.
\end{remark}

\begin{remark}
Every Jacob's ladder
\bdis
\vp_1(T)=\frac 12\vp(t)
\edis
where $\vp_t$ is the exact solution of the nonlinear integral equation (\ref{1.5}) is the asymptotic solutions of the following nonlinear \emph{integro-iterations} equation
\be \label{2.5}
\begin{split}
& \frac 1U\int_T^{T+U}\prod_{k=0}^n\left|\zeta\left(\frac 12+ix^k(t)\right)\right|^2{\rm d}t= \\
& = \prod_{k=0}^n\frac{1}{x^k(T+U)-x^k(T)}\int_{x^k(T)}^{x^k(T+U)}\left|\zf\right|^2{\rm d}t
\end{split}
\ee
of the new kind, where (comp. (\ref{2.1}))
\bdis
x^0(t)=t,\ x^1(t)=x(t),\ x^2(t)=x(x(t)),\ \dots \ ,
\edis
i.e. the function $x^k(t)$ is the $k$-th iteration of the function $x(t),\ t\geq T_0[\vp_1]$ (comp. (\ref{2.5}) with \cite{3}, (11.1), (11.4), (11.6), (11.8) and \cite{4}, (2.6))).
\end{remark}

\section{Proof of Theorem}

\subsection{}

We use the formulae (see \cite{2}, (6.2))

\be \label{3.1}
t-\vp_1(t)\sim (1-c)\pi(t)\sim (1-c)\frac{t}{\ln t},\quad t\to\infty .
\ee

\begin{remark}
The fundamental geometric property of the set of Jacob's ladders is expressed by the formula (\ref{3.1}). Namely, the difference of the abscissa and the ordinate of the point $[t,\vp_1(t)]$ of the curve
$y=\vp_1(t)$ is asymptotically equal to $(1-c)\pi(t)$.
\end{remark}

We have (see (\ref{3.1}))
\be \label{3.2}
\begin{split}
& \vp_1^1(t)-\vp_1^2(t)\sim (1-c)\frac{\vp_1^1(t)}{\ln\vp_1^1(t)}, \\
& \vp_1^2(t)-\vp_1^3(t)\sim (1-c)\frac{\vp_1^2(t)}{\ln\vp_1^2(t)} , \\
& \vdots \\
& \vp_1^n(t)-\vp_1^{n+1}(t)\sim (1-c)\frac{\vp_1^n(t)}{\ln\vp_1^n(t)},\ t\in [T,T+U]
\end{split}
\ee
for arbitrary fixed $n\in\mbb{N}$ and (see (\ref{3.1}), (\ref{3.2}))

\be \label{3.3}
\begin{split}
& t\sim \vp_1^1(t)\sim\vp_1^2(t)\sim \dots \sim \vp_1^{n+1}(t),\ t\to\infty , \\
& t>\vp_1^1(t)>\vp_1^2(t)>\dots >\vp_1^{n+1}(t) \ .
\end{split}
\ee

Next we have (see (\ref{3.2}), (\ref{3.3}))

\be \label{3.4}
\begin{split}
& \vp_1^1(t)-\vp_1^2(t)\sim (1-c)\frac{t}{\ln t} , \\
& \vp_1^2(t)-\vp_1^3(t)\sim (1-c)\frac{t}{\ln t}, \\
& \vdots \\
& \vp_1^n(t)-\vp_1^{n+1}(t)\sim (1-c)\frac{t}{\ln t} ,
\end{split}
\ee
and, consequently, by the addition of (\ref{3.1}) and (\ref{3.4}), we obtain

\be \label{3.5}
\begin{split}
& t-\vp_1^{n+1}(t)\sim (1-c)(n+1)\frac{t}{\ln t} , \\
& \vp_1^{n+1}(t)\sim t\left\{ 1-\frac{(1-c)(n+1)}{\ln t}\right\},\ 0<1-c<1 , \\
& \vp_1^{n+1}(t)>\left(1-\frac{\epsilon}{2}\right)t \left\{ 1-\frac{(1-c)(n+1)}{\ln t}\right\}>(1-\epsilon)t\geq (1-\epsilon)T .
\end{split}
\ee
(see (\ref{2.2})), i.e. from (\ref{3.5}) by (\ref{3.3}) (the second line) the formula (\ref{2.4}) follows. Especially, the following holds true (see (\ref{3.3}), (\ref{3.5}))

\be \label{3.6}
(1-\epsilon)T<\vp_1^{n+1}(T)<T .
\ee

\subsection{}

Let (see \cite{3}, (9.1), (9.2))

\be \label{3.7}
\tilde{Z}^2(t)=\frac{{\rm d}\vp_1(t)}{{\rm d}t},\ \vp_1(t)=\frac 12\vp(t),\ t\geq T_0[\vp_1]
\ee
where

\be \label{3.8}
\begin{split}
& \tilde{Z}^2(t)=\frac{Z^2(t)}{2\Phi^\prime_\vp[\vp(t)]}=\frac{\left|\zf\right|^2}{\left\{ 1+\mcal{O}\left(\frac{\ln\ln t}{\ln t}\right)\right\}\ln t} , \\
& t\in [T,T+U],\ U\in \left(\left. 0,\frac{T}{\ln T}\right]\right. \ .
\end{split}
\ee

If we use the formula (\ref{3.7}) for the iterations (\ref{2.1}) we obtain

\be \label{3.9}
\prod_{k=0}^n \tilde{Z}^2[\vp_1^k(t)]=\frac{{\rm d}\vp_1^1}{{\rm d}t}\frac{{\rm d}\vp_1^2}{{\rm d}\vp_1^1}\cdots \frac{{\rm d}\vp_1^n}{{\rm d}\vp_1^{n-1}}\frac{{\rm d}\vp_1^{n+1}}{{\rm d}\vp_1^n}=
\frac{{\rm d}\vp_1^{n+1}}{{\rm d}t}
\ee
by the rule for differentiation of a composite function. Next, by the integration of (\ref{3.8}) we obtain

\be \label{3.10}
\int_T^{T+U}\prod_{k=0}^n\tilde{Z}^2[\vp_1^k(t)]{\rm d}t=\vp_1^{n+1}(T+U)-\vp_1^{n+1}(T) .
\ee

\subsection{}

It follows easily from (\ref{3.7}) that

\bdis
\begin{split}
& \int_{\vp_1^k(T)}^{\vp_1^k(T+U)}\tilde{Z}^2(t){\rm d}t=\vp_1(\vp_1^k(T+U))-\vp_1(\vp_1^k(T))= \\
& = \vp_1^{k+1}(T+U)-\vp_1^{k+1}(T),\quad k=0,1,\dots , n ,
\end{split}
\edis
i.e.

\bdis
\frac{1}{\vp_1^k(T+U)-\vp_1^k(T)} \int_{\vp_1^k(T)}^{\vp_1^k(T+U)}\tilde{Z}^2(t){\rm d}t=\frac{\vp_1^{k+1}(T+U)-\vp_1^{k+1}(T)}{\vp_1^k(T+U)-\vp_1^k(T)} ,
\edis
and, consequently,

\be \label{3.11}
\prod_{k=0}^n\frac{1}{\vp_1^k(T+U)-\vp_1^k(T)}\int_{\vp_1^k(T)}^{\vp_1^k(T+U)}\tilde{Z}^2(t){\rm d}t=\frac{\vp_1^{n+1}(T+U)-\vp_1^{n+1}(T)}{U} .
\ee
Hence, from (\ref{3.10}), (\ref{3.11}) we obtain the \emph{exact} formula

\be\label{3.12}
\begin{split}
& \frac 1U\int_T^{T+U}\prod_{k=0}^n \tilde{Z}^2[\vp_1^k(t)]{\rm d}t = \\
& = \prod_{k=0}^n\frac{1}{\vp_1^k(T+U)-\vp_1^k(T)}\int_{\vp_1^k(T)}^{\vp_1^k(T+U)}\tilde{Z}^2(t){\rm d}t .
\end{split}
\ee

\subsection{}

First of all, we have (see (\ref{2.2}), (\ref{3.1}))

\bdis
\begin{split}
& \left( 1-\frac{1}{2n+1}\right)(1-c)\frac{T}{\ln T}<T+U-\vp_1^1(T+U)<\left( 1+\frac{1}{2n+1}\right)(1-c)\frac{T}{\ln T} , \\
& \left( 1-\frac{1}{2n+1}\right)(1-c)\frac{T}{\ln T}<T-\vp_1^1(T+U)<\left( 1+\frac{1}{2n+1}\right)(1-c)\frac{T}{\ln T} .
\end{split}
\edis

Next, we have ($0<1-c<1$)

\bdis
\begin{split}
& |\{T+U-\vp_1^1(T+U)\}-\{T-\vp_1^1(T)\}|< \\
& < \left( 1+\frac{1}{2n+1}\right)(1-c)\frac{T}{\ln T}-\left( 1-\frac{1}{2n+1}\right)(1-c)\frac{T}{\ln T} = \\
& =\frac{2}{2n+1}(1-c)\frac{T}{\ln T}<\frac{2}{2n+1}\frac{T}{\ln T} ,
\end{split}
\edis
i.e.

\be \label{4.1}
\vp_1^1(T+U)-\vp_1^1(T)-U<\frac{2}{2n+1}\frac{T}{\ln T} ,
\ee
and (see (\ref{2.2}))

\bdis
\begin{split}
& 0<\vp_1^1(T+U)-\vp_1^1(T)<\frac{2}{2n+1}\frac{T}{\ln T}+U\leq \frac{2}{2n+1}\frac{T}{\ln T}+\frac{T}{\ln^2T}< \\
& < \frac{2}{2n+1}\frac{T}{\ln T}+\frac{1}{2n+1}\frac{T}{\ln T}=\frac{3}{2n+1}\frac{T}{\ln T} .
\end{split}
\edis
Hence (see (\ref{2.2}))

\be \label{4.2}
U\leq \frac{T}{\ln^2T} \ \Rightarrow \ \vp_1^1(T+U)-\vp_1^1(T)<\frac{3}{2n+1}\frac{T}{\ln T} .
\ee
Similarly, from the formula (see (\ref{3.4}))

\bdis
\vp_1^1(t)-\vp_1^2(t)\sim (1-c)\frac{t}{\ln t},\ t\to\infty
\edis
we obtain (comp. (\ref{4.1})) that (see (\ref{4.2}))

\bdis
\vp_1^2(T+U)-\vp_1^2(T)<\frac{2}{2n+1}\frac{T}{\ln T}+\vp_1^1(T+U)-\vp_1^1(T)<\frac{5}{2n+1}\frac{T}{\ln T} .
\edis

Next, if the estimate (the function $\vp_1^k$ is increasing)

\bdis
(0<)\ \vp_1^k(T+U)-\vp_1^k(T)<\frac{2k+1}{2n+1}\frac{T}{\ln T}
\edis
holds true then we obtain, by a similar way, that

\bdis
\vp_1^{k+1}(T+U)-\vp_1^{k+1}(T)<\frac{2}{2n+1}\frac{T}{\ln T}+\frac{2k+1}{2n+1}\frac{T}{\ln T}<\frac{2(k+1)+1}{2n+1}\frac{T}{\ln T} .
\edis

Hence, the following estimate holds true

\be \label{4.3}
\vp_1^k(T+U)-\vp_1^k(T)<\frac{2k+1}{2n+1}\frac{T}{\ln T}\leq \frac{T}{\ln T},\ k=1,\dots , n.
\ee
(Comp. the condition for $U$ in (\ref{3.8})).

\subsection{}

Let us remind that by (\ref{3.6}) we have the following inequalities

\bdis
(1-\epsilon)T<\vp_1^{n+1}(T)<T+U,\quad U\in \left(\left. 0,\frac{T}{\ln^2T}\right.\right] ,
\edis
i.e. for every
\bdis
T'\in (\vp_1^{n+1}(T),T+U)
\edis
the following holds true

\be \label{3.13}
\ln T'=\ln T+\mcal{O}(1),\quad T\to\infty .
\ee

Now, if we use the mean-value theorem in (\ref{3.12}) we obtain (see (\ref{3.3}), (\ref{3.8}), (\ref{4.3}), (\ref{3.13}))

\be \label{3.14}
\begin{split}
& \int_T^{T+U}\prod_{k=0}^n\tilde{Z}^2[\vp_1^k(t)]{\rm d}t\sim
\frac{1}{\ln^{n+1}T}\int_T^{T+U}\prod_{k=0}^n\left|\zeta\left(\frac 12+i\vp_1^k(t)\right)\right|^2{\rm d}t ,
\end{split}
\ee
\be \label{3.15}
\begin{split}
& \prod_{k=0}^n\frac{1}{\vp_1^k(T+U)-\vp_1^k(T)}\int_{\vp_1^k(T)}^{\vp_1^k(T+U)}\tilde{Z}^2(t){\rm d}t\sim \\
& \sim \frac{1}{\ln^{n+1}T}\prod_{k=0}^n\frac{1}{\vp_1^k(T+U)-\vp_1^k(T)}\int_{\vp_1^k(T)}^{\vp_1^k(T+U)}\left|\zf\right|^2{\rm d}t .
\end{split}
\ee
Hence, from (\ref{3.12}) by (\ref{3.14}), (\ref{3.15}) the asymptotic formula (\ref{2.3}) follows.

\thanks{I would like to thank Michal Demetrian for helping me with the electronic version of this work.}

\end{document}